\documentclass[english]{ourlematema}
\pdfoutput=1
\usepackage{amsmath, amsthm, amssymb, hyperref, color}
\usepackage{MnSymbol} 
\usepackage{graphicx}
\usepackage{caption}
\usepackage[all]{xypic}
\usepackage{verbatim}
\usepackage{chemfig, chemnum}
\usepackage{tikz}
\usepackage[linesnumbered,lined,commentsnumbered]{algorithm2e}
\usepackage{caption}
\usepackage{subcaption}
\usepackage{float}
\usepackage{makecell}
\usepackage{array}
\usepackage{fancyvrb}
\usepackage{enumitem}
\usepackage{standalone}

\newcommand{\bR}{\mathbb R}
\newcommand{\bC}{\mathbb C}

\newcommand{\bP}{\mathbb P}

\newcommand{\cC}{\mathcal C}

\newcommand{\cO}{\mathcal O}

\newcommand{\cV}{\mathcal V}
\newcommand{\cW}{\mathcal W}

\renewcommand{\phi}{\varphi}

\newcommand{\Gr}{\text{Gr}}

\newcommand{\codim}{\text{codim}}

 \title{Chow-Lam Recovery}
 
  \author{Elizabeth Pratt}
  \address{%
  University of Berkeley California, Berkeley \\
\email{epratt@berkeley.edu}
}

  \author{Kristian Ranestad}
  \address{%
  University of Oslo, Oslo \\
\email{ranestad@math.uio.no}
}

 \date{2024/9/2}

\begin{document}
\maketitle
\begin{abstract}
\noindent 
We study the conditions under which a subvariety of the Grassmannian may be recovered from certain of its linear projections.  In the special case that our Grassmannian is projective space, this is equivalent to asking when a variety can be recovered from its Chow form; the answer is ``always" by work of Chow in 1937 \cite{CW}.
In the general Grassmannian setting, the analogous question is when a variety can be recovered from its Chow-Lam form. We give both necessary conditions for recovery and families of examples where, in contrast with the projective case, recovery is not possible. 
\end{abstract}

\section{Introduction}

The Grassmannian $\Gr(k,n)$ parameterizes $(k-1)$-dimensional subspaces of $\bP^{n-1}.$ We are interested in projections between Grassmannians induced by projections of the underlying projective spaces. That is, fix an $n \times r$ matrix $Z$ of rank $r$ and consider the map
\begin{align}\label{eq:zproj}
\begin{split}
    Z: \Gr(k, n) & \dashrightarrow \Gr(k,r) \\
    [M] & \mapsto [MZ]
    \end{split}
\end{align}
where $M$ is a $k \times n$ matrix and $[M]$ is its rowspan. Geometrically, one may think of projecting each $(k-1)$-plane from $\bP^{n-1}$ to $\bP^{r-1}.$ Indeed, when $k=1$ we recover linear projections between projective spaces. 

In this paper, $\cV$ will be a subvariety of $\Gr(k,n)$ of dimension $k(r-k) -1$ for some $r \leq n.$ Let $Z(\cV)$ denote the closure of the image of $\cV$ under the projection map induced by $Z$. Then for a general $Z$ the variety $Z(\cV)$ is expected to be a hypersurface in the target Grassmannian, and is thus cut out by a single equation. 

We ask the question: can  $\cV$ be recovered from the data of \emph{all} possible $Z(\cV)$? Let $Z^{-1}(Z(\cV))$ denote the pre-image of $Z(\cV)$ in $\Gr(k,n).$ In general, we have the containment
\begin{equation}\label{eq:zcap}
    \cV  \ \subseteq \ \bigcap_{Z} Z^{-1}(Z(\cV)).
\end{equation}
If $k=1$, equality holds in \ref{eq:zcap}. In general the right side may be strictly larger.

Such linear projections arise in the construction of the amplituhedron in particle physics. The context is that of particle scattering for $N=4$ super-symmetric Yang-Mills. Recall that the Grassmannian $\Gr(k,n)$ has a \emph{Pl\"ucker embedding} into projective space $\bP^{\binom{n}{k}-1}$, given by taking the $k \times k$ minors of each $k \times n$ matrix representative of a subspace. The projective coordinates on $\bP^{\binom{n}{k}-1}$ are called the \emph{Pl\"ucker coordinates} of $\Gr(k,n).$ We define the \emph{positive Grassmannian} to be the semi-algebraic subset of $\Gr_{\bR}(k,n)$ where the Pl\"ucker coordinates are non-negative. When $Z$ has strictly positive maximal minors, any $Z$-projection of the positive Grassmannian into $\Gr(k,r)$ is called the \emph{amplituhedron} (the images for different such $Z$ are combinatorially equivalent). The amplituhedron is used in calculations of scattering amplitudes for gluons \cite{AHT}. 

The boundaries of the positive Grassmannian are known as \emph{positroid varieties} \cite{Williams}. Thomas Lam studied positroid varieties and their $Z$-projections in \cite{Lam}, and his work motivated the later definition of the Chow-Lam form in \cite{PS}. We will return to positroid varieties in Section \ref{sec:schubert} of this paper.

We now define the Chow-Lam form as in \cite{PS} and explain the connection to projections. The input data is a variety $\cV \subset \Gr(k,n)$ of dimension $k(r-k)-1$ for some $r \leq n.$ For a subspace $Q \subset \bP^{n-1},$ let $\Gr(k, \hat Q)$ denote the Grassmannian of $(k-1)$-spaces contained in $Q.$ The \emph{Chow-Lam locus} is the variety
\begin{equation*}
    \mathcal{CL}_\cV := \{Q \in \Gr(n-r+k, n): \cV \cap \Gr(k, \hat Q) \neq \varnothing \}.
\end{equation*}
In words, the Chow-Lam locus consists of linear spaces $Q$ which contain as a subspace some $(k-1)$-space $P,$ where $P$ is a point in $\cV.$ By Lemma 3.2 of \cite{PS}, $\mathcal{CL}_\cV$ is a proper irreducible subvariety of $\Gr(n-r+k, n).$ When $\mathcal{CL}_\cV$ is a hypersurface, it is defined by a single equation, which we call the \emph{Chow-Lam form} and denote $CL_\cV$. If $\mathcal{CL}_\cV$ is not a hypersurface, we set $CL_\cV := 1.$

\begin{exa}[Curve in $\Gr(2,4)$]
    Let $\cV$ be a curve in the Grassmannian $\Gr(2,4).$  Here $k = 2, n = 4, $ and $r = 3.$ Thus the Chow-Lam locus lives in $\Gr(3,4) = (\bP^3)^\vee.$ It consists of planes containing a line in $\bP^3$, such that the line is a point of $\cV.$ 

    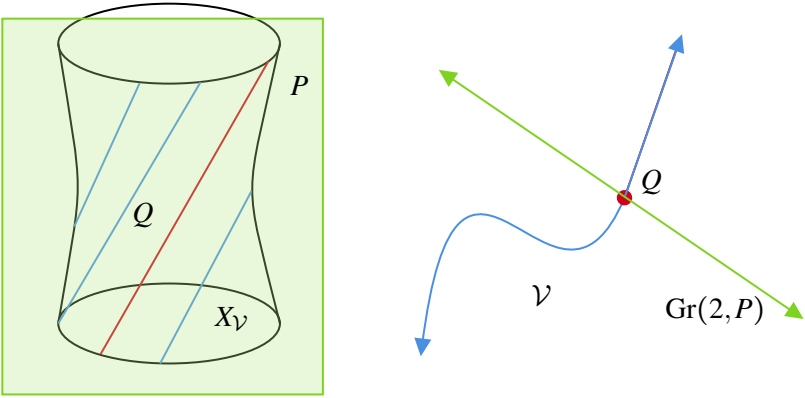
\begin{figure}[!h]
		\begin{minipage}{0.5\textwidth}
			\centering
			\tikzset{every picture/.style={line width=0.75pt}} 

\begin{tikzpicture}[x=0.75pt,y=0.75pt,yscale=-1,xscale=1]

\draw   (249,63) .. controls (249,51.95) and (273.85,43) .. (304.5,43) .. controls (335.15,43) and (360,51.95) .. (360,63) .. controls (360,74.05) and (335.15,83) .. (304.5,83) .. controls (273.85,83) and (249,74.05) .. (249,63) -- cycle ;
\draw   (249,203) .. controls (249,191.95) and (273.85,183) .. (304.5,183) .. controls (335.15,183) and (360,191.95) .. (360,203) .. controls (360,214.05) and (335.15,223) .. (304.5,223) .. controls (273.85,223) and (249,214.05) .. (249,203) -- cycle ;
\draw    (249,63) .. controls (262,144.5) and (262,125.5) .. (249,203) ;
\draw    (360,63) .. controls (342,145.5) and (341,126.5) .. (360,203) ;
\draw [color={rgb, 255:red, 74; green, 144; blue, 226 }  ,draw opacity=1 ]   (320.25,82.25) -- (249,203) ;
\draw [color={rgb, 255:red, 208; green, 2; blue, 27 }  ,draw opacity=1 ]   (354.25,71.75) -- (270.25,218.25) ;
\draw [color={rgb, 255:red, 74; green, 144; blue, 226 }  ,draw opacity=1 ]   (345.75,136.75) -- (300,223) ;
\draw [color={rgb, 255:red, 74; green, 144; blue, 226 }  ,draw opacity=1 ]   (289.75,82.75) -- (257,154.5) ;
\draw  [color={rgb, 255:red, 126; green, 211; blue, 33 }  ,draw opacity=1 ][fill={rgb, 255:red, 184; green, 233; blue, 134 }  ,fill opacity=0.3 ] (221,50.5) -- (381.5,50.5) -- (381.5,238.5) -- (221,238.5) -- cycle ;

\draw (363.83,77.83) node [anchor=north west][inner sep=0.75pt]   [align=left] {$\displaystyle P$};
\draw (284.67,141) node [anchor=north west][inner sep=0.75pt]   [align=left] {$\displaystyle Q$};
\draw (326,193.17) node [anchor=north west][inner sep=0.75pt]   [align=left] {$\displaystyle X_{\mathcal{V}}$};

\end{tikzpicture}
		\end{minipage}
		\begin{minipage}{0.45\textwidth}
			\centering
			\tikzset{every picture/.style={line width=0.75pt}} 

\begin{tikzpicture}[x=0.75pt,y=0.75pt,yscale=-1,xscale=1]

\draw [color={rgb, 255:red, 74; green, 144; blue, 226 }  ,draw opacity=1 ]   (209.55,216.79) .. controls (227.98,69.82) and (276.35,224.67) .. (311,142) ;
\draw [shift={(209,221.33)}, rotate = 276.64] [fill={rgb, 255:red, 74; green, 144; blue, 226 }  ,fill opacity=1 ][line width=0.08]  [draw opacity=0] (8.93,-4.29) -- (0,0) -- (8.93,4.29) -- cycle    ;
\draw [color={rgb, 255:red, 208; green, 2; blue, 27 }  ,draw opacity=1 ]   (338.6,63) -- (311,142) ;
\draw [shift={(311,142)}, rotate = 109.26] [color={rgb, 255:red, 208; green, 2; blue, 27 }  ,draw opacity=1 ][fill={rgb, 255:red, 208; green, 2; blue, 27 }  ,fill opacity=1 ][line width=0.75]      (0, 0) circle [x radius= 3.35, y radius= 3.35]   ;
\draw [color={rgb, 255:red, 74; green, 144; blue, 226 }  ,draw opacity=1 ]   (338.8,63.03) -- (311,142) ;
\draw [shift={(339.8,60.2)}, rotate = 109.4] [fill={rgb, 255:red, 74; green, 144; blue, 226 }  ,fill opacity=1 ][line width=0.08]  [draw opacity=0] (8.93,-4.29) -- (0,0) -- (8.93,4.29) -- cycle    ;
\draw [color={rgb, 255:red, 126; green, 211; blue, 33 }  ,draw opacity=1 ]   (220.48,79.69) -- (398.19,200.98) ;
\draw [shift={(400.67,202.67)}, rotate = 214.31] [fill={rgb, 255:red, 126; green, 211; blue, 33 }  ,fill opacity=1 ][line width=0.08]  [draw opacity=0] (8.93,-4.29) -- (0,0) -- (8.93,4.29) -- cycle    ;
\draw [shift={(218,78)}, rotate = 34.31] [fill={rgb, 255:red, 126; green, 211; blue, 33 }  ,fill opacity=1 ][line width=0.08]  [draw opacity=0] (8.93,-4.29) -- (0,0) -- (8.93,4.29) -- cycle    ;

\draw (317.67,126.33) node [anchor=north west][inner sep=0.75pt]   [align=left] {$\displaystyle Q$};
\draw (330,189.67) node [anchor=north west][inner sep=0.75pt]   [align=left] {$\displaystyle {\rm Gr}(2,P)$};
\draw (263,184.33) node [anchor=north west][inner sep=0.75pt]   [align=left] {$\displaystyle \mathcal{V}$};

\end{tikzpicture}
		\end{minipage}
		\caption{Geometry in $\bP^3$ (left) and $\Gr(2,4)$ (right)}
		\label{fig:twoschubert}
	\end{figure}
    
    Let $X_\cV$ be the surface in $\bP^3$ swept out by all of the lines in $\cV$. Then any tangent plane to $X_\cV$ contains all lines through the point of tangency; thus the dual surface $X_\cV^\vee$ is contained in $\mathcal{CL}_\cV.$ In fact, since they are irreducible varieties of the same dimension, we have that the Chow-Lam locus is precisely $X_\cV^\vee.$
\end{exa}

When $k=1,$ we recover the definition of the \emph{Chow form} from classical algebraic geometry. This was pioneered by Cayley in 1860 for curves in three-space \cite{Cayley} and generalized by Chow and van der Waerden in 1937 \cite{CW}.

If $A, B$ are subvarieties of $\bP^n,$ we define their \emph{join} $A \vee B$ to be the closure of the union of all lines $\overline{ab}$ spanned by distinct points $a$ in $A$ and $b$ in $B.$ The fibers of the projections $Z$ of the form \ref{eq:zproj} can be described using join.  Let $K_Z\in \Gr(n-r,n)$ be the projectivized kernel of the $r\times n$ matrix $Z^T$.

\begin{lemma} \label{lem:inverseimg} 
    Fix a linear space $P \in \Gr(k,n),$ and assume that $K_Z$ does not intersect $P$. Then $P \vee K_Z$ is in $\mathcal{CL}_\cV$ if and only if $P$ is in $Z^{-1}(Z(\cV)).$ 
\end{lemma}

\begin{proof}
   Observe that $Z^{-1}(Z(P))$, as a subset of $\Gr(k,n)$, is exactly the Grassmannian $\Gr(k,\widehat{P\vee K_Z})$. So $Z(P)$ is in $Z(\cV)$ only if $P \vee K_Z$ contains as a subspace some $Q \in \cV.$  
\end{proof}

Inspired by Lemma \ref{lem:inverseimg}, we define for a given $(k(r-k)-1)$-dimensional variety $\cV\subset \Gr(k,n)$, the algebraic set
\begin{equation}
    W_\cV := 
    \{ P \in \Gr(k,n) \ : \  Q \in \mathcal{CL}_\cV \text{ whenever } Q\supset P\} 
\end{equation}
and call it the \emph{recovered variety} of $\cV.$ (We will call it a variety, but in general it may be reducible.) We remind the reader that the dimension of $Q$ is $n-r+k-1.$ Note that if $Q\supset P$, then $Q\supset P\vee A$ for every $A \in \Gr(n-r,\hat Q)$, and each such $A$ is the projectivized kernel $K_Z$ for some $n \times r$ matrix $Z$. This gives us the following corollary.

\begin{cor}
    The varieties $W_\cV$ and $\bigcap_{Z} Z^{-1}(Z(\cV))$ are equal.
\end{cor}

In particular, $W_\cV$ contains $\cV.$ Our question about projections becomes: are there points in the recovered variety which are not in $\cV$? 
 
For the Chow form, we have that $W_\cV = \cV$. That is, any projective variety can be uniquely recovered from its Chow form.  In fact, the Chow form was originally studied as a way to represent a projective variety by a single equation in Pl\"ucker coordinates \cite{CW}. However, this is not true for the Chow-Lam form: two different varieties may have the same Chow-Lam form, as in Example \ref{eg:ruledquadric}.

\begin{exa}[A ruled quadric]\label{eg:ruledquadric}
    Let $X$ be a quadric surface in $\bP^3.$ Then $X$ has two rulings, which give two curves $\cV$ and $\cV'$ in $\Gr(2,4).$ The Chow-Lam locus of each of these curves is the dual variety to $X.$ 
\end{exa}


The structure of this paper is as follows. In Section \ref{sec:computation} we will review the theory of Chow-Lam forms and their computation. In Section \ref{sec:criteria} we will discuss a general criterion for a point in $\Gr(k,n)$ to lie in $W_\cV.$ The remaining three sections are devoted to examples. Section \ref{sec:curves and surfaces} discusses curves and surfaces, Section \ref{sec:ruled} describes families of varieties with the same Chow-Lam form, and Section \ref{sec:schubert} describes $W_\cV$ when $\cV$ is a Schubert hyperplane section of the Grassmannian.

\section{Chow-Lam Forms and Computation}\label{sec:computation}
To compute Chow-Lam forms and their recovery, it is convenient to use two types of Pl\"ucker coordinates.  An $(m-1)$-subspace $P$ of $\bP^{n-1}$ can be represented as the rowspan of a $m \times n$ matrix, or as the projectivized kernel of an $(n-m) \times n$ matrix.  In the first case, the $m \times m$ minors of this matrix are called the \emph{dual Pl\"ucker coordinates} of $P$ and denoted $q_I(P)$, where $I \subset \binom{[n]}{m}.$ In the latter case, the maximal minors are called the \emph{primal Pl\"ucker coordinates} and denoted $p_I(P),$ where $I \subset \binom{[n]}{n-m}.$

If $m$ is close to $n$ then it is more convenient to use primal coordinates. Primal and dual Pl\"ucker coordinates are related up to sign by complementary indices. For example, the ten Pl\"ucker coordinates on $\Gr(3, 5)$ are

\begin{equation}
\label{eq:primaldual35}
 \begin{matrix}
 p_{12} & p_{13} & p_{14} & p_{15} & p_{23} & p_{24} & p_{25} & p_{34} &
p_{35} & p_{45}, \\ 
q_{345} & -q_{245} & q_{235} & -q_{234} & q_{145} & -q_{135} & q_{134} & q_{125} & -q_{124} & q_{123}.
\end{matrix}
\end{equation}

Suppose that $\cV \subset \Gr(k,n)$ is a variety of dimension $k(r-k)-1$ for some $r \geq k.$ Then the Chow-Lam form of $\cV$ can be computed by first computing the equations for the incidence variety $$\Phi = \{(P, Q) : P \subset Q\} \subset \Gr(k,n) \times \Gr(n-r+k, n).$$ We then intersect $\Phi$ with $\cV \times \Gr(n-r+k, n)$ and project to the latter factor. Algebraically, this corresponds to taking the sum of the ideal of $\Phi$ with the ideal of $\cV$ and eliminating the variables in $\Gr(k,n).$ 

\begin{exa}[Equations of incidence variety]\label{eg:inc}
Suppose $k = 2, n = 4,$ and $r = 3.$ The incidence variety $\Phi$ will consist of pairs of lines and planes in $\bP^3$ such that the plane contains the line. If $p_i$ are the primal coordinates of the plane and $q_{ij}$ are the dual coordinate of the line, then the ideal of $\Phi$ is generated by
\begin{align}
    \begin{split}
        q_{12}p_1-q_{23}p_3-q_{24}p_4, \quad q_{12}p_2+q_{13}p_3+q_{14}p_4, \\
    q_{13}p_1+q_{23}p_2-q_{34}p_4, \quad q_{14}p_1+q_{24}p_2+q_{34}p_3.
    \end{split}
\end{align}
\end{exa}

Given the Chow-Lam form $CL_\cV$ of a variety $\cV$, the variety $W_\cV$ can be computed as follows. We expand 
the dual Pl\"ucker coordinates of $A \vee P\in \Gr(n-r+k,n)$ in terms of dual Pl\"ucker coordinates of $A\in \Gr(n-r,n)$ and $P\in\Gr(k,n),$ then collect coefficients of monomials in the $q_I(A).$ When $A \vee P\in \mathcal{CL}_\cV$, these coefficients are polynomials in the variables $q_I(P)$ which define the variety $W_\cV.$ Note that this process produces not only a variety, but a scheme; indeed, one may define the \emph{recovered scheme} $\mathcal{W}_\cV$ as the scheme defined by these polynomials. In the classical case of the Chow form for a variety in projective space, cf. \cite[Theorem 1.14]{Catanese}, the scheme $\mathcal{W}_\cV$ is equal to $\cV$ exactly at the smooth points of $\cV$ and has embedded components at the singular points. We can see this in Example \ref{eg:quintic}. In the Chow-Lam case the situation is more complicated: even for a smooth subvariety $\cV$ of the Grassmannian with $W_\cV = \cV,$ there may be embedded components in $\cW_\cV$, as in Example \ref{eg:scheme}.

\begin{exa}[Singular quintic in $\bP^3$]\label{eg:quintic}
		Consider the curve $\cV \subset \bP^3$ given parametrically by $[s:t] \mapsto [s^5:s^4t:s^3t^2:t^5].$ Let $\bP^3$ have coordinates $x_1,x_2,x_3,x_4.$ Then $\cV$ has a singularity of Jacobian degree $4$ at $[0:0:0:1]$ and is otherwise smooth; that is, the vanishing of the Jacobian and the equation for $\cV$ define an ideal of length $4$ whose variety is that point. Because the curve is given to us parametrically, we can also compute the Chow form parametrically: it is the closure of rowspans of matrices
        \begin{align*}
            L_s = \begin{bmatrix}
                1 & t & t^2 & t^5 \\
                a & b & c & d
            \end{bmatrix}
        \end{align*}
        and obtained by eliminating $t, a, b, c, d.$ In primal Pl\"ucker coordinates this gives
        \begin{align*}
        \begin{split}
            C_\cV = p_{14}^5-& p_{13}^3p_{14}p_{24}+3p_{12}p_{13}p_{14}^2p_{24}+p_{12}^3p_{24}^2\\
            & + p_{13}^4p_{34}-p_{12}^3p_{23}p_{34}-4p_{12}p_{13}^2p_{14}p_{34}+2p_{12}^2p_{14}^2p_{34}.
        \end{split}
        \end{align*}
		We view the $p_{ij}$ as primal Pl\"ucker coordinates of the line $x \vee a$ where $x$ and $a$ are points in $\bP^3.$ Then $x$ is in $\cV$ if and only if $x \vee a$ is in $\cC_\cV$ for all $a \in \bP^3.$ Thus we can make the change of coordinates $p_{12} = x_3a_4-x_4a_3,$ etc.  Collecting coefficients of monomials in the $a$ variables gives us an ideal which is the intersection of two primary components. One is the ideal of $\cV,$ and the other is the degree $16$ ideal 
		\[(x_4^2, x_3^3x_4, x_2x_3^2x_4, x_2^2x_3x_4, x_2^3x_4, x_3^4, x_2x_3^3, x_2^2x_3^2, x_2^3x_3-x_1x_3^2x_4, x_2^4-x_1x_3^3)\]
		whose variety is the point $[0:0:0:1].$
	\end{exa}

\begin{exa}[Curve in $\Gr(2,4)$]\label{eg:scheme}
		Consider the curve $\cV$ in $\Gr(2,4)$ given by tangent lines to the quintic in Example \ref{eg:quintic}. In dual Pl\"ucker coordinates, this curve is given by
		\begin{align*}
			\cV&  = V(q_{24}^2-4q_{23}q_{34}, 3q_{14}q_{24}-10q_{13}q_{34}, q_{13}q_{24}-6q_{12}q_{34}, 8q_{14}^2-25q_{23}q_{24}, \\
			& q_{23}q_{14}-5q_{12}q_{34}, 15q_{23}^2-4q_{13}q_{14}, 2q_{13}q_{23}-3q_{12}q_{24}, 5q_{13}^2-9q_{12}q_{14}).
		\end{align*}
		Then $\cV$ has a singular locus consisting of two points $[0:0:0:0:0:1]$ and $[1:0:0:0:0:0],$ each of Jacobian degree $2.$ The Chow-Lam form of $\cV$ in $\Gr(3,4)=\bP^3$ is the degree $7$ form given in primal Pl\"ucker coordinates by
  \begin{align*}
      \begin{split}
          CL_{\cV} = 16p_2^3p_3^4& +108p_1^2p_3^5-128p_2^4p_3^2p_4\\
     & -900p_1^2p_2p_3^3p_4+256p_2^5p_4^2+2000p_1^2p_2^2p_3p_4^2+3125p_1^4p_4^3.
      \end{split}
  \end{align*}
		
		Here $CL_{\cV}$ can be computed by using the equations in Example \ref{eg:inc} and eliminating. The ideal of the recovered scheme has three primary components of degrees $14,42,$ and $52.$ Their corresponding varieties are the curve $\cV$ and the two singular points. This is an example where $W_{\cV}$ and $\cV$ are the same as algebraic sets, but scheme-theoretically different even at the smooth points, as seen by the degree.
	\end{exa}

To keep our discussion simple, 
we will restrict ourselves to considering the recovered variety $W_\cV$ as an algebraic set. We end this section by noting that the dimension of $W_\cV$ may be higher than that of $\cV.$

\begin{exa}[High dimension $W_\cV$]
Choose a basis $e_1, ..., e_n$ of $\bC^n.$ Let $\cW \subset \Gr(k,n)$ be the Schubert variety of subspaces meeting $e_n,$ which we now view as a point in $\bP^{n-1}$. Then $\cW$ is isomorphic to the Grassmannian $\Gr(k-1, n-1).$ Let $\cV$ be a generic $(k-1)$-dimensional subvariety of $\cW.$ We claim that $\cW$ and $\cV$ have the same $Z$-projections for a generic matrix $Z$ of dimension $n \times (k+1).$ Thus $\cW$ will appear as a component in $W_\cV.$

To see this, let $Z_n$ be the last row of $Z.$ Then we may view $Z_n$ as a point in $\bP^k.$ The projection $Z(\cW)$ consists of all $(k-1)$-spaces in $\bP^k$ meeting the point $Z_n.$ So $Z(\cW)$ is a hyperplane in $\Gr(k, k+1)=(\bP^k)^\vee.$   To conclude it suffices to show that $Z(\cV)$ has the same dimension as $\cV$ for generic $\cV \subset \cW.$ But this follows by Proposition 4.8 of \cite{Lam}, which gives cohomological conditions on $\cV$ for the dimension of the projection to agree with the dimension of the variety.

How much bigger is $\cW$? It has codimension $n-k$ and dimension $(k-1)(n-k)$ in $\Gr(k,n).$ Recall that $\cV$ has dimension $k+1.$ Thus as $k$ increases, we get a family of examples which show that $W_\cV$ can be arbitrarily large compared to $\cV.$
\end{exa}

\section{A Recovery Criterion} \label{sec:criteria}
Fix a variety $\cV \subset \Gr(k,n)$ of dimension $k(r-k)-1$ for some $r.$ In this section, we will establish a general criterion for a point $P \in \Gr(k,n)$ to be in $W_\cV.$ 

Fix a linear space $P \subset \bP^{n-1}$ of dimension $k-1.$ Then we may stratify $\cV$ as follows. Let 
\[\cV_{P,i} = \{P' \in \cV : \dim P' \cap P \geq i\} \subset \Gr(k,n) \]
be the Schubert variety of spaces in $\cV$ which meet $P$ in dimension at least $i.$ We take the empty set to have dimension $-1.$ Then we have 
\[\cV = \cV_{P,-1} \supset \cV_{P,0} \supset ... \supset \cV_{P,k-1},\]
where the last term is nonempty exactly when $P$ itself is in $\cV.$  We begin with a lemma on the dimension of certain Schubert varieties, then state our main theorem.
Let $Q \subset \bP^{n-1}$ be a subspace of dimension $m$ and assume $l>m$.  Then we denote by $\Omega(Q)\subset \Gr (l,n)$ the Schubert variety of subspaces that contain $Q$.

\begin{lemma}[Dimension]\label{lem:dim}
    Suppose $Q \subset \bP^{n-1}$ is a subspace of dimension $m$. Then the variety $\Omega(Q)$ has dimension $(l-m-1)(n-l)$ in $\Gr(l,n).$  
\end{lemma}
\begin{proof}  $\Omega(Q)$ is isomorphic to the Grassmannian $\Gr(l-m-1,n-m-1)$, so has dimension $(l-m-1)(n-l)$ in $\Gr(l,n).$  
    \end{proof}

\begin{thm}[Recovery criterion]\label{thm:resid}
  Let $\cV\subset\Gr(k,n)$ be a subvariety of dimension $k(r-k)-1$.    If a linear space $P\in\Gr(k,n)$ is in the recovery $W_\cV$, then for some $0 \leq i \leq k-1,$ the variety $\cV_{P,i}$ has dimension at least $(k-i-1)(r-k).$
\end{thm}

\begin{proof} 
    The proof is essentially a dimension argument.  Consider the Schubert variety $\Omega(P) \subset \Gr(n-r+k,n)$ of $(n-r+k-1)$-spaces containing $P.$ Observe that a linear space $P$ is in the recovery if and only if every point of $\Omega(P)$ is also a point of $\Omega(P')$ for some $P' \in \cV.$ 
    
    So we consider the incidence
    \[
    I_P=\{(R,P') \ | \ R\in \Omega(P)\cap \Omega(P'), P'\in \cV\} \ \subset \ \Gr(n-r+k,n)\times \Gr(k,n),
    \]
    with projections $\pi:I_P\to \Omega(P)$ and $\pi':I_P\to \cV.$ Then $P$ is in the recovery if and only if $\pi$ is surjective.The latter is the case only if the dimension of $I_P$ is at least $(n-r)(r-k)={\rm dim }\Omega(P).$
   Now $I_P$ may have several components since the fibers of $\pi'$ may differ in dimensions.
   The fiber $(\pi')^{-1}(P')$ is the variety $\Omega(P\wedge P')\subset \Gr(n-r+k,n)$.
   The dimension of this fiber depends on the dimension of the intersection $P\cap P'$.  Observe that
   $$\dim P'\vee P = \dim P' + \dim P - \dim (P' \cap P).$$
   When $P'\in \cV_{P,i}\setminus \cV_{P,i+1}$, the span $P\wedge P'$ of $P$ and $P'$ has dimension $2k-i-2$, so $\Omega(P\cap P')$ has dimension $(n-r-k+i+1)(r-k)$.  Therefore the preimage over $\cV_{P,i}$ dominates $\Omega(P)$ only if
   $\dim\cV_{P,i}+ (n-r-k+i+1)(r-k)\geq (n-r)(r-k),$ i.e. when $\dim\cV_{P,i}\geq (k-i-1)(r-k)$.


\end{proof}

For example, let $\cV$ be a surface in $\Gr(3,n),$ whose points represent planes in $\bP^{n-1}.$ Then $\dim \cV = 2 = 3(r-3) - 1,$ so $r = 4.$ If a plane $P$ is in $W_\cV,$ then by Theorem \ref{thm:resid} at least one of the following must be true:
\begin{itemize}
    \item $(i = 0):$ $P$ meets each $P' \in \cV$ in a point, or
    \item $(i = 1):$ $P$ meets a one-dimensional family of $P' \in \cV$ in a line, or
    \item $(i = 2):$ $P$ is in $\cV$, i.e. meets a point in $\cV$ in a plane.
\end{itemize}
The dimension condition may be satisfied for multiple $i$. For example, suppose all planes $P'$ in $\cV$ contain a fixed line $L$ in $\bP^{n-1}.$ Then for any plane $P$ not in $\cV$ containing $L,$ the dimension condition is satisfied for both $i=0$ and $i = 1.$

Theorem \ref{thm:resid} gives a necessary condition for $P$ to be in the recovery, but it is not sufficient in general. Even if the variety $\cV_{P, i}$ has high dimension, the map $\pi$ may have high-dimensional fibers and so may not be surjective.

\begin{exa}[Counterexample to the converse of Theorem \ref{thm:resid}]
    This example will show that the condition in Theorem \ref{thm:resid} is not sufficient for $P$ to be in $W_\cV.$ Let $n$ be large, and let $e_i$ be the $i$th standard basis vector. Fix the plane $Q := e_1 \vee e_2 \vee e_3$ and the $3$-space $P := e_1 \vee e_2 \vee e_3 \vee e_4.$ Next, we let $\cV_1$ be a  $4$-dimensional subvariety of $\Gr(4,n)$ with the property that all $P'$ in $\cV_1$
    \begin{enumerate}[label=\roman*)]
    \item meet $Q$ in at least a line
    \item meet $e_{n-1} \vee e_n$ in at least a point.
    \end{enumerate}
    Finally, choose a general $7$-dimensional variety $\cV\subset \Gr(4,n)$ containing $\cV_1.$ Then $r = 6$ and the Chow-Lam locus $\mathcal{CL}_\cV$ lives in $\Gr(n-2,n).$ The dimension of $\cV_{P,1}$ is at least $(k-1-1)(r-k) = 4,$ which per Theorem \ref{thm:resid} is necessary for $P$ to be in the recovery; indeed, the dimension of $\cV_{P,0}$ is lower than $6$ if $\cV$ is chosen generically.  However, we will see that $P$ is not in the recovery.

    To see this, take the codimension two subspace $R \subset \bP^{n-1}$ given by the span of all vectors except $e_{n-1}$ and $e_n.$ Then $R$ contains $P$, but by construction cannot contain any $P'$ in $\cV_1.$  If $\cV$ is chosen generically to contain $\cV_1,$ then $R$ will not contain any $P'$ in $\cV$ either. Thus  $\pi:I_P\to \Omega(P)$ is not surjective, and $P$ is not in the reecovery $W_\cV.$
\end{exa}

\section{Chow-Lam Locus as a Dual Variety}\label{sec:curves and surfaces}

Given a variety $\cV \subset \Gr(k,n),$ we define $X_\cV \subset \bP^{n-1}$ to be the union of all $(k-1)$-spaces in $\cV,$ or equivalently the image under the projection to $\bP^{n-1}$ of the incidence 
\[\{(x,P) \ | \ x\in P\subset \bP^{n-1}, P\in \cV\} \ \subset \ \bP^{n-1}\times \Gr(k,n).
\]
The following proposition is a generalization of Example \ref{eg:ruledquadric} of the quadric surface.

\begin{prop}\label{prop:dual}
    Let $\cV$ be a subvariety of $\Gr(k,n)$ of dimension $k-1.$  Then the dual $X_\cV^\vee$ is contained in $\mathcal{CL}_\cV,$ with equality whenever they have the same dimension.
    In this case a $(k-1)$-space $P$ belongs to the recovery $W_\cV$ if and only if its dual linear space $P^\bot$ is contained in $X_\cV^\vee$.
\end{prop}

\begin{proof}
    Note that if $n$ is less than $2k,$ then $X_\cV = \bP^{n-1}$ and its dual is empty, meaning the proposition has no content. For the rest of the proof we suppose that $n \geq 2k$. Here $r=k+1,$ so the Chow-Lam locus is a hypersurface in $\Gr(n-1,n) \cong (\bP^{n-1})^\vee.$ Take any tangent space to $X_\cV$ at a smooth point $x.$ Then it must contain all linear spaces in $X_\cV$ through $x.$ In particular it contains all $(k-1)$-dimensional linear spaces through $x$ that belong to $\cV$.  But by definition, $x$ is in at least one $(k-1)$-space in $\cV.$ Thus any hyperplane that contains such a tangent space, i.e. any hyperplane that belong to  $X_\cV^\vee$, is in the Chow-Lam locus. Both  $X_\cV^\vee$ and $\mathcal{CL}_\cV$ are irreducible, so we have equality whenever their dimensions coincide. In this case the last statement of the proposition is immediate.
\end{proof}



Let us apply Theorem \ref{thm:resid} to curves. When a curve $\cV\subset \Gr(k,n)$ has a Chow-Lam form, then $k(r-k)-1=1$. Thus either $k=1$ and we are in the ordinary case of Chow forms, or $k=2$ and $r=3$ and the Chow-Lam form defines a hypersurface in $\Gr(n-1,n)$. In the latter case, the Proposition \ref{prop:dual} applies, and the Chow-Lam form defines the dual variety $X_\cV^\vee$ whenever the latter is non-degenerate. 

More generally, when $\cV\subset \Gr(k,n)$ has a Chow-Lam form and $k(r-k)-1=p-1$, where $p$ is prime, then $k=1$ and we are in the ordinary case of Chow forms, or $k=p,r=p+1$ and the Chow-Lam locus lives in $\Gr(n-1,n)$.
So we have the following corollary of Proposition \ref{prop:dual}:

\begin{cor} Fix a variety $\cV\subset \Gr(k,n), k>1$ of dimension $p-1$ for $p$ prime. Then $\cV$ has a Chow-Lam form only if $k = p$ and $\mathcal{CL}_\cV \subset (\bP^n)^\vee.$ 
In this case, a $(p-1)$-space $P$ is in the recovery $W_\cV$ if and only if $P^{\bot}\subset X_\cV^\vee$.
\end{cor}
In the case $p=2$, i.e. when $\cV$ is a curve with a Chow-Lam form, then $k=2$. With this assumption we can prove a stronger statement.
\begin{cor} Let $\cV\subset \Gr(2,n)$ be a curve such that $X_\cV$ is not a cone.
Then a line $L\in\bP^{n-1}$ is in the recovery 
$W_\cV$ if and only if $L\subset X_\cV$.
\end{cor}
\begin{proof} When the surface scroll $X_\cV$ is not a cone, the dual variety $X_\cV^\vee$ is a hypersurface which coincides with the Chow-Lam locus. If $L$ is a line contained in the surface $X_\cV$, then every hyperplane that contains $L$ will contain a line that belongs to $\cV$, so $L$ will belong to the recovery $W_\cV$.  For the converse, we will show that if $L$ is not contained in $X_\cV$ then we can find a hyperplane that contains $L$ but no line in $\cV$. 

First note that since $X_\cV$ is not a cone, only finitely many lines in $\cV$ intersects $L$.  Now, consider a general codimension two linear space $Q$ that contain $L$. It will intersect $X_\cV$ in finitely many points.  So there are only finitely many lines in $\cV$ that intersect $Q$. But then there are only finitely many hyperplanes that contain both $Q$ and a line in $\cV$. Therefore a general hyperplane containing $Q$ will contain no line in $\cV.$


\end{proof}





\begin{exa}[Hirzebruch surface]
    Consider the Hirzebruch surface $$X = \bP(\cO_{\bP^1}(1) \oplus \cO_{\bP^1}(a)).$$ Since $X$ is the projectivization of a bundle, there is a map $\pi: X \to \bP^1$ whose fibers are isomorphic to $\bP^1$. Then $X$ is ruled by these fibers, all of which meet the curve $L\subset X$ corresponding to the line subbundle $\cO_{\bP^1}(1)$. Embed $X$ into $\bP^{a+2},$ and let $\cV$ be the curve in $\Gr(2,a+3)$ whose points are the lines in the ruling of $X.$ Then the image of $L$ is a line and the recovery of $\cV$ will contain this line as a point in $\Gr(2,n)\setminus \cV.$ 

    For a concrete example, let $a = 1$. Then $X$ is the closure of the embedding
    \begin{align*}
        (\bC^*)^2 & \to \bP^4 \\
        (s:t) & \mapsto [1:s:s^2:st:s^2t].
    \end{align*}
    Each line in the ruling is given by fixing $s.$ Parametrically, they are rowspans of $$L_s = \begin{bmatrix}
        1 & s & s^2 & 0 & 0 \\
        1 & s & s^2 & s & s^2
    \end{bmatrix}.$$ These lines form a curve $\cV \subset \Gr(2,5)$ given by the Pl\"ucker relations and 
    $$V(q_{45}, q_{34}-q_{25}, q_{24}-q_{15}, q_{23}, q_{13}, q_{12}).$$ The Chow-Lam form may be calculated parametrically by adding two rows to $L_s,$ and is $p_3p_4^2-p_2p_4p_5+p_1p_5^2.$ 
    
    The recovered scheme will be the original curve as well as a fat point of length $25,$ whose radical ideal consists of all Pl\"ucker coordinates except for $q_{45}.$ This is the line
    spanned by $[0:0:0:1:0]$ and $[0:0:0:0:1]$ obtained in the closure of the embedding by allowing $t$ to approach infinity.
    
\end{exa}

In a surface scroll $X_\cV$, 
there may even be two lines $L_1$ and $L_2$  in the scroll that are not rulings.  In this case every line in the ruling is in the linear span of $L_1$ and $L_2$, i.e. a $\bP^3$, so  $n=4$ and the scroll $X_\cV$ is a surface of degree $d_1+d_2$ in $\bP^3$ with multiplicity $d_1$ along $L_1$ and multiplicity $d_2$ along $L_2$.  In fact, we explain the multiplicities as follows. We have that $\cV\subset \Gr(2,4)$ is contained in the quadric surface of lines meeting $L_1$ and $L_2$.  The surface $X_\cV$ is the birational image of $\bP_\cV(\mathcal{L}_1+\mathcal{L}_2)$ in $\bP^3$. Here $\mathcal{L}_1$ is a line bundle of degree $d_i$ and $L_i$ is the image of $\cV$ by a basepoint free pencil of sections of $\mathcal{L}_i$. The degree of $X_\cV$ is then computed  on $\bP_\cV(\mathcal{L}_1+\mathcal{L}_2)$, or in $\bP^3$. 
The map of $\cV$ to $L_i$ has degree $d_i$, so $X_\cV$ must have multiplicity $d_i$ along $L_i$. The degree of $X_\cV$ in $\bP^3$ is computed as a sum of intersection multiplicities with a line.  A general line that intersect $L_1$ and $L_2$ will intersect $X_\cV$ only in these lines, and the sum of intersection multiplicities is the sum of the multiplicities of $X_\cV$ along these lines, i.e. $d_1+d_2$.  

\begin{exa}[Maps to a curve and a line]\label{eg:notjoin}
 Consider a curve $\cC$ with a $d_1:1$ map $\phi_1$ to $\bP^1$ and a birational map $\phi_2$ to $\bP^{n-3}$ with image a curve $C$ of degree $d_2.$ Next, we embed $\bP^{n-3}$ and $\bP^1$ as non-intersecting linear spaces in $\bP^{n-1}$. Let $X$ be obtained by joining, for each $c\in \cC$, the points $\phi_1(c)$ and $\phi_2(c)$ in $\bP^{n-1}.$ Then the two maps $\phi_1$ and $\phi_2$ define a map $\cC\to \Gr(2,n)$ with image $\cV$ such that $X=X_{\cV}$.  The recovery of $\cV$ contains the line $\phi_1(C)$ in $X$ that does not belong to $\cV$. 

 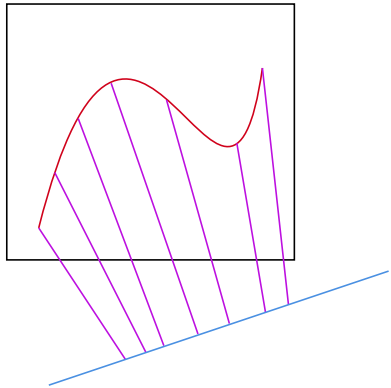
\begin{figure}[!h]
 \centering
\scalebox{0.8}{\tikzset{every picture/.style={line width=0.75pt}} 

\begin{tikzpicture}[x=0.75pt,y=0.75pt,yscale=-1,xscale=1]

\draw  [line width=0.75]  (100,20) -- (280,20) -- (280,180) -- (100,180) -- cycle ;
\draw [color={rgb, 255:red, 208; green, 2; blue, 27 }  ,draw opacity=1 ][line width=0.75]    (120,160) .. controls (178,-68.5) and (238,215.5) .. (260,60) ;
\draw [color={rgb, 255:red, 189; green, 16; blue, 224 }  ,draw opacity=1 ][line width=0.75]    (200,80) -- (239.33,219.83) ;
\draw [color={rgb, 255:red, 189; green, 16; blue, 224 }  ,draw opacity=1 ][line width=0.75]    (260,60) -- (276.33,208.17) ;
\draw [color={rgb, 255:red, 189; green, 16; blue, 224 }  ,draw opacity=1 ][line width=0.75]    (120,160) -- (174.33,242.33) ;
\draw [color={rgb, 255:red, 189; green, 16; blue, 224 }  ,draw opacity=1 ][line width=0.75]    (144.71,91.57) -- (198.33,233.67) ;
\draw [color={rgb, 255:red, 189; green, 16; blue, 224 }  ,draw opacity=1 ][line width=0.75]    (165.33,69.5) -- (219.67,226.33) ;
\draw [color={rgb, 255:red, 189; green, 16; blue, 224 }  ,draw opacity=1 ][line width=0.75]    (244,107) -- (262,213.17) ;
\draw [color={rgb, 255:red, 189; green, 16; blue, 224 }  ,draw opacity=1 ][line width=0.75]    (130,125.33) -- (187,237.67) ;
\draw [color={rgb, 255:red, 74; green, 144; blue, 226 }  ,draw opacity=1 ][line width=0.75]    (339,187) -- (126.33,258.33) ;

\end{tikzpicture}
 }
 \caption{Construction of a rational scroll}
     
 \end{figure}
    
\end{exa}
Analogously, we find surfaces  $\cV\subset \Gr(3,n), n\geq 10$ for which $W_\cV\neq \cV$.

\begin{exa}(Maps to a surface and a plane).
Suppose an abstract (not embedded) surface $\mathcal{S}$ has two very ample linear systems $L_1,L_2$ and a basepoint free $2$-dimensional linear system $L_3$. Let $X_\cV$ be the $4$-dimensional $\bP^2$-scroll spanned by the images of $\mathcal{S}$ by the sum of the three linear systems (composed possibly with a general linear projection to $\bP^{n-1}$ with $n\geq 10$). This gives us an embedding $\mathcal{S}\to\cV\subset \Gr(3,n)$ such that the plane image of $\mathcal{S}$ by $L_3$ belongs to the $W_\cV$ but not to $\cV$. 
 \end{exa}
\section{Multi-ruled Varieties}\label{sec:ruled}
	In this section we generalize Example \ref{eg:ruledquadric} and construct subvarieties of $\Gr(k,n)$ which have the same Chow-Lam form. We begin with an immediate corollary of Proposition \ref{prop:dual}.

\begin{cor}[Multi-ruled varieties]\label{cor:multiruled}
    Let $\cV_1, \cV_2$ be subvarieties of $\Gr(k,n)$ of dimension $k-1.$ Let $X_i \subset \bP^{n-1}$ be the closure of the union of linear spaces in $\cV_i$ for $i = 1,2.$ Suppose that the duals $X_i^\vee$ are non-degenerate. Then $X_1 = X_2$ if and only if $\cV_1$ and $ \cV_2$ have the same Chow-Lam locus.
    \end{cor}


    \begin{exa}[Segre varieties]
    Consider the Segre embedding $\Sigma_K := \bP^{k-1} \times \bP^{k-1} \hookrightarrow \bP^{k^2-1}.$ Let $\cV_1$ and $\cV_2$ be the two $(k-1)$-dimensional subvarieties  of $\Gr(k, k^2)$ that parameterize $(k-1)$-spaces in $\Sigma_K$. That is, $\cV_1$ parameterizes the $(k-1)$-spaces in the second factor, 
    \begin{align*}
        \cV_1 = \biggl\{ \text{Seg}(\{p\} \times \bP^{k-1}) \in \rm{Gr}(k,k^2) \ : \ p \in \bP^{k-1}\biggr\},
    \end{align*}
    while $\cV_2$ parameterizes the $(k-1)$-spaces in the first factor. Then the varieties $X_i$ as in the Corollary \ref{cor:multiruled} coincide with $\Sigma_K$. 
    By Theorem 0.1 in \cite{GKZ}, the dual variety of $\Sigma_K$ is non-degenerate. Thus Corollary \ref{cor:multiruled} tells us that $\cV_1$ and $\cV_2$ both have Chow-Lam locus equal to $\Sigma_k^\vee.$
    \end{exa}
    
    \begin{exa}[Multi-Segre varieties]
    We can also consider Segre embeddings with more factors than two: for example, $\bP^1 \times \bP^1 \times \bP^1 \times \bP^1 \hookrightarrow \bP^{15}.$ Let $\cV_i$ be the variety in $\Gr(2,16)$ corresponding to embedding the $i$th factor, and letting points on the other factors vary. In other words, $\cV_i$ is a family of lines parameterized by $\bP^1 \times \bP^1 \times \bP^1.$ Then $\dim \cV_i = 3,$ and so $\cV_i$ has a Chow-Lam locus which lives in $\Gr(16-4+3,16) = \Gr(15,16).$ It is exactly the dual variety to the Segre variety. 

    More generally, for any $k \geq 2$ and $i \geq 0$, the Segre embedding $(\bP^{k-1})^{ki+2}$ is ruled by $ki+2$ families of $(k-1)$-spaces. This gives an example of $ki+2$ varieties in $\Gr(k, N)$ with the same Chow-Lam locus, where $N = k^{ki+2}-1.$ Each variety has dimension  $(ki+1)(k-1) = k(r -k)-1,$ where $r = ki + 1 -i + k.$ By Theorem 0.1 of \cite{GKZ} the dual variety of the Segre embedding is non-degenerate, and thus the Chow-Lam form of each $\cV_i$ is exactly this dual variety.
    \end{exa}

    \begin{exa}[Linear determinantal hypersurfaces]\label{eg:detsurfaces}

 The Segre variety defines two varieties isomorphic to $\bP^{k-1}$ in $\Gr(k,k^2)$.  Their Chow-Lam locus lives in $\Gr(k^2-1,k^2) = (\bP^{k^2-1})^\vee$ and is the dual hypersurface to the Segre variety. We can view $\bP^{k^2-1}$ as the space of $k\times k$-matrices in which the Segre variety is the variety of matrices of rank one.  The dual space may then also be interpreted as $(k\times k)$-matrices, under the nondegenerate bilinear map $$(A,B)\mapsto Trace(A\cdot B).$$
 If $A$ is a matrix of rank one, then a general tangent hyperplane to the Segre variety at $A$ is a matrix $B$ of corank one such that the matrix product $A\cdot B$ vanishes.  So the dual hypersurface to the Segre variety is the matrices of corank one, i.e. the determinantal hypersurface $S_k \subset \bP^{k^2-1}$ of degree $k$. We call a linear combination of rows in the matrix defining $S_k$ a ``generalized row." It will be a vector of linear forms. Each generalized row in the matrix defining $S_k$ vanishes on a linear space of codimension $k$ contained in $S_k$, and similarly for generalized columns. Thus $S_k$ contains two $\bP^{k-1}$-families of codimension $k$ subspaces.
The variety $S_k$ has a degenerate dual variety, namely the Segre variety $\Sigma_{k}$ itself.  

We can get additional examples from linear projections of $\Sigma_k.$ The dual variety of a linear section  $L \cap S_k$ of $S_k$ coincides with the projection of $\Sigma_{k}$ from $L^{\bot}$, as long as $L^{\bot}$ does not intersect $\Sigma_{k}$, which is possible only if $L$ has dimension at least $2k-1$.  In fact, when $L^{\bot}$ does not intersect the Segre variety, the tangent hyperplanes to the projected variety are precisely the tangent hyperplanes in $\bP^{k^2-1}$ that contain $L^{\bot}$. These hyperplanes correspond to the points on $S_k$ that lie in $L$.   

Suppose $L$ has dimension $m=2k-1$.The projection of the Segre variety gives rise to two varieties isomorphic to $\bP^{k-1}$ in $\Gr(k,2k)$ with Chow-Lam form the determinantal hypersurface $D_k:=L\cap S_k$.  Since $m=2k-1$ both the projection of the Segre variety and 
$D_k$ are  hypersurfaces.  They are dual hypersurfaces that both  have two rulings of $(k-1)$-spaces.  In fact, the two families of codimension $k$ subspaces in $S_k$ define, after intersection with $L$, two families of $(k-1)$-spaces on $D_k=L\cap S_k$.  So this case gives rise to two pairs of varieties in $\Gr(k,2k)$ with coinciding Chow-Lam form.  For the two varieties defined by rows and columns of the $(k\times k)$-matrices restricted to $L$, the Chow-Lam form defines the projected Segre variety and has degree $\binom{2k-2}{k-1}$, the degree of the Segre variety.  For the two varieties defined by the two factors in the Segre variety, the Chow-Lam form defines the determinantal hypersurface $D_k$.

 \end{exa}
 \begin{exa}[Cubic surface]
     This is a special case of Example \ref{eg:detsurfaces} in which $k=3$ and thus $m = 5$. Consider a $3 \times 3$ matrix $M$ filled with linear forms in six variables. Let $S$ denote the determinantal hypersurface cut out by $\det M$. Then $S$ is ruled by two families $\cV_1, \cV_2$ of planes. A plane $[a:b:c]$ in the first family is given by setting the linear combination $aR_1 + bR_2 + cR_3 $ of the rows to zero. Similarly, a plane in the second family is given by setting a linear combination of the columns to zero. These give two subvarieties of $\Gr(3,6).$ 

     The Chow-Lam locus is the dual variety $S^\vee.$ We may think of $S^\vee$ as a projection of the Segre variety $\Sigma_3$ from $(\bP^8)^\vee$ to $(\bP^5)^\vee.$ The Chow-Lam form is the defining equation of $S^\vee,$ and has degree $\binom{4}{2} = 6.$ 

     There is a connection to cubic surfaces as follows. We may calculate the degree of $S^\vee$ by intersecting it with a line in $(\bP^5)^\vee$, or equivalently by intersecting $S$ with a generic pencil of hyperplanes in $\bP^5.$ The base locus $L$ of this pencil is a $\bP^3.$ The points on $S^\vee$ are hyperplanes in $\bP^5$ that contains a plane in $S$. So the intersection points of $S^\vee$ with $L^{\bot}$ correspond exactly to the planes in $S$ which meet $L$ in a line. 
     
     Now, the restriction $S|_L$ will be a cubic surface in $L = \bP^3$.
     So, one can calculate the degree of $S^\vee$ using this description, as follows. Write the $3 \times 3$ matrix $M$ as a $3 \times 3 \times 6$ tensor. Write $L$ as the rowspan of a $4 \times 6$ matrix, which we call $M_L.$ Then we may restrict to $L$ by multiplying our tensor by $M_L,$ giving a $3 \times 3 \times 4$ tensor $T$.  For a concrete example, we might have that the slices of $T$ are
     \begin{align*}
         T_x = \begin{bmatrix}
         0 & 1 & 1 \\ 
         2 & -1 & -3 \\
         0 & 3 & 3
         \end{bmatrix} \ , \  T_y = 
         \begin{bmatrix}
         2 & -3 & -1 \\ -1 & -3 & 0 \\ 3 & 0& -2
         \end{bmatrix} \\
         T_z = \begin{bmatrix}
             0 & 0 & 2 \\ 
             -2 & 0 & 1\\ 
             -1 & -2 & 3
         \end{bmatrix} \ , \ T_w = \begin{bmatrix}
             3 & -2 & 0 \\ 0 & 3 & -3 \\ 2 & 1 & -1
         \end{bmatrix}.
     \end{align*}
    Then the cubic surface $S|_L$ is given by the vanishing of
    \begin{align*}
       F :=  \det \begin{bmatrix}
            2y+3w & x-3y-2w & x-y+2z \\
            2x-y-2z & -x-3y+3w & -3x+z-3w \\ 3y-z+2w & 3x-2z+w & 3x-2y+3z-w
        \end{bmatrix}.
    \end{align*}
    The intersection points of $S^{\vee}$ with $L^\bot$ are the exactly set of points $[a:b:c]$ in $\bP^2$ such that $\{aR_1 + bR_2 + cR_3 = 0\}$ corresponds to a line rather than a point. In other words, these are the points where the matrix $\begin{bmatrix}
        a & b & c
    \end{bmatrix} \cdot T$ drops rank. In our example, we have
    \begin{align*}
        \begin{bmatrix}
        a & b & c
    \end{bmatrix} \cdot T =  \begin{bmatrix}
        b+c & 2a-3b-c & 2c & 3a-2b \\
        2a-b-3c & -a-3b & -2a+c& 3b-3c \\
        3b+3c&  3a-2c & -a-2b+3c& 2a+b-c
    \end{bmatrix}.
    \end{align*}

    By Porteous' formula, the degree of the degeneracy locus of $\begin{bmatrix}
        a & b & c
    \end{bmatrix} \cdot T$ is $6,$ which is our desired degree.  For example, in our case the degeneracy locus consists of six complex points cut out by
    \begin{small}
    \begin{align*}153669b^3-9842a^2c+126290abc+211487b^2c+29792ac^2-152462bc^2-21896c^3, \\ 153669ab^2+127798a^2c+223799abc+47621b^2c-79510ac^2+109111bc^2-87722c^3, \\ 153669a^2b+151069a^2c-331777abc+513941b^2c-112573ac^2-104261bc^2-104150c^3, \\ 153669a^3-600557a^2c+298853abc+542006b^2c+314438ac^2-752480bc^2+231640c^3.
    \end{align*}
    \end{small}

    If we had multiplied $T$ with $[a \ b \ c]$ on a different side, we would have a matrix whose degeneracy locus consists of those generalized columns which define lines rather than points. Thus we obtain two families of lines, six coming from the rows and six from the columns. These form a ``double six" on the cubic surface defined by $F$: namely, they do not pairwise intersect and any line from one family meets exactly five lines from another family. 
\end{exa}
    \section{Schubert Arrangements} \label{sec:schubert}
    In this section we explore in detail the case where $\cV$ is a linear section of $\Gr(k,n)$ by Schubert divisors. We will see that in this case the variety $W_{\cV}$ contains additional linear components not in $\cV$. 
 
    Fix $n, k$ with $k \leq n,$ and a partition $\lambda=(\lambda_0,...,\lambda_{k-1})$ which fits in a $k \times (n-k)$ box. Let $E_0 \subset ... \subset E_{n-1}$  be the standard flag, where $E_i$ is the linear space spanned by the first $i+1$ coordinate points. We define the \emph{Schubert variety $\Omega_{\lambda}$} to be
    \[\Omega_{\lambda} = \{L \in \Gr(k,n):{\rm dim}\; L \cap E_{n-k+i-\lambda_i} \geq i \text{ for all } i \}.\]

    For convenience we may sometimes write our partitions in the form $m^{a_m}\ldots 1^{a_1}$ where $a_j$ is the number of parts of size $j.$ For example, $\lambda = (3,3,2,2,2,2,0)$ would be written as $3^22^4.$ For more on Schubert varieties, see \cite{Fulton}.

    Let $H$ be a linear space of codimension $k$ in $\bP^{n-1}.$ Let $\Omega_1(H) \subset \Gr(k,n)$ be the variety of subspaces in $\bP^{n-1}$ which meet $H$. In our notation, $\Omega_1(H)$ is isomorphic to the Schubert variety $\Omega_1,$ i.e. the Schubert variety of codimension $1$ in $\Gr(k,n)$. In this section, we will consider Schubert varieties of the form $\Omega_1^j = \Omega_1(H_1) \cap ... \cap \Omega_1(H_j),$ for $j$ general linear spaces $H_i$, with $j=k(n-r)+1$ for some $k<r<n$. We will refer to the Schubert varieties in $W_{\cV} \setminus \cV$ as \emph{recovered Schubert components}. Note that our notation for $\Omega_1^j$ does not depend on $n$; we will show that given $j,$ there is a set of recovered Schubert components which will appear for any $n$ large enough.

    \begin{exa}[Threefold in $\Gr(2,5)$]
	Choose $H_1, H_2, H_3$ general planes in $\bP^4.$ Consider the threefold $\cV = \Omega_1^3$ of lines in $\bP^4$ meeting all three planes. The Chow-Lam locus of $\cV$ will live in $\Gr(3,5).$ Then $W_\cV$ is the variety of all lines $L$ such that for all planes $Q$ containing $L,$ $Q$ also contains a common transversal to $H_1, H_2$, and $H_3.$ We claim that $W_\cV$ also contains three extra components of each of two equivalence classes:
	\begin{enumerate}[label=\roman*)]
		\item The Schubert variety $\Omega_3$  of lines meeting $H_i \cap H_j$
		\item The Schubert variety $\Omega_{2^2}$ of lines contained in $H_i$.
	\end{enumerate}

    Let $Q \cong \bP^2$ be any plane. For a generic choice of $Q,$ the planes $H_1, H_2, H_3$ will intersect $Q$ in three points $p_1, p_2, p_3,$ which do not have a common transversal. To show that $L$ is in the recovery, we need to check that the condition $Q \supset L$ forces $p_1, p_2, p_3$ to have a common transversal in $Q.$

	For case (i), suppose $L$ meets $H_1 \cap H_2$ and consider the geometry in a generic plane $Q \cong \bP^2$ containing $L.$ Then $Q$ also contains the point $p = H_1 \cap H_2 \in L.$ Thus the line $\overline{pp_3}$ is a common transversal to $H_1, H_2, H_3.$ 
	
	For case (ii), suppose $L$ is contained in $H_1$ and consider a generic plane $Q$ containing $L.$ Then $H_1$ intersects $Q$ not in a point, but a line: the line $L$ itself. Then the line $\overline{p_2p_3}$ will intersect $L$ in $Q$ by B\'ezout's theorem. Thus $\overline{p_2p_3}$ will be a common transversal to $H_1, H_2, H_3.$ Both these cases are illustrated in Figure \ref{fig:twoschubert}.
	
	\begin{figure}[!h]
		\begin{minipage}{0.48\textwidth}
			\centering
			\tikzset{every picture/.style={line width=0.75pt}} 

\begin{tikzpicture}[x=0.75pt,y=0.75pt,yscale=-1,xscale=1]
				
				\draw   (180,40) -- (400,40) -- (400,240) -- (180,240) -- cycle ;
				\draw    (280,100) -- (320,180) ;
				\draw [shift={(320,180)}, rotate = 63.43] [color={rgb, 255:red, 0; green, 0; blue, 0 }  ][fill={rgb, 255:red, 0; green, 0; blue, 0 }  ][line width=0.75]      (0, 0) circle [x radius= 3.35, y radius= 3.35]   ;
				\draw [shift={(280,100)}, rotate = 63.43] [color={rgb, 255:red, 0; green, 0; blue, 0 }  ][fill={rgb, 255:red, 0; green, 0; blue, 0 }  ][line width=0.75]      (0, 0) circle [x radius= 3.35, y radius= 3.35]   ;
				\draw    (260.89,61.79) -- (280,100) ;
				\draw [shift={(260,60)}, rotate = 63.43] [color={rgb, 255:red, 0; green, 0; blue, 0 }  ][line width=0.75]    (10.93,-3.29) .. controls (6.95,-1.4) and (3.31,-0.3) .. (0,0) .. controls (3.31,0.3) and (6.95,1.4) .. (10.93,3.29)   ;
				\draw    (320,180) -- (339.11,218.21) ;
				\draw [shift={(340,220)}, rotate = 243.43] [color={rgb, 255:red, 0; green, 0; blue, 0 }  ][line width=0.75]    (10.93,-3.29) .. controls (6.95,-1.4) and (3.31,-0.3) .. (0,0) .. controls (3.31,0.3) and (6.95,1.4) .. (10.93,3.29)   ;
				\draw    (201.94,119.51) -- (358.06,80.49) ;
				\draw [shift={(360,80)}, rotate = 165.96] [color={rgb, 255:red, 0; green, 0; blue, 0 }  ][line width=0.75]    (10.93,-3.29) .. controls (6.95,-1.4) and (3.31,-0.3) .. (0,0) .. controls (3.31,0.3) and (6.95,1.4) .. (10.93,3.29)   ;
				\draw [shift={(200,120)}, rotate = 345.96] [color={rgb, 255:red, 0; green, 0; blue, 0 }  ][line width=0.75]    (10.93,-3.29) .. controls (6.95,-1.4) and (3.31,-0.3) .. (0,0) .. controls (3.31,0.3) and (6.95,1.4) .. (10.93,3.29)   ;
				
				\draw (275,62) node [anchor=north west][inner sep=0.75pt]   [align=left] {$\displaystyle  p = p_1= p_2$};
				\draw (330,160) node [anchor=north west][inner sep=0.75pt]   [align=left] {$\displaystyle p_{3}$};
				\draw (201,82) node [anchor=north west][inner sep=0.75pt]   [align=left] {$\displaystyle L$};

			\end{tikzpicture}

   
   \caption*{$L$ is in $\Omega_3$}
		\end{minipage}
  		\begin{minipage}{0.48\textwidth}
			\centering
						\tikzset{every picture/.style={line width=0.75pt}} 
			
			\begin{tikzpicture}[x=0.75pt,y=0.75pt,yscale=-1,xscale=1]
				
				\draw   (200,40) -- (420,40) -- (420,240) -- (200,240) -- cycle ;
				\draw    (320,140) -- (300,180) ;
				\draw [shift={(300,180)}, rotate = 116.57] [color={rgb, 255:red, 0; green, 0; blue, 0 }  ][fill={rgb, 255:red, 0; green, 0; blue, 0 }  ][line width=0.75]      (0, 0) circle [x radius= 3.35, y radius= 3.35]   ;
				\draw [shift={(320,140)}, rotate = 116.57] [color={rgb, 255:red, 0; green, 0; blue, 0 }  ][fill={rgb, 255:red, 0; green, 0; blue, 0 }  ][line width=0.75]      (0, 0) circle [x radius= 3.35, y radius= 3.35]   ;
				\draw    (359.11,61.79) -- (340,100) ;
				\draw [shift={(360,60)}, rotate = 116.57] [color={rgb, 255:red, 0; green, 0; blue, 0 }  ][line width=0.75]    (10.93,-3.29) .. controls (6.95,-1.4) and (3.31,-0.3) .. (0,0) .. controls (3.31,0.3) and (6.95,1.4) .. (10.93,3.29)   ;
				\draw    (300,180) -- (280.89,218.21) ;
				\draw [shift={(280,220)}, rotate = 296.57] [color={rgb, 255:red, 0; green, 0; blue, 0 }  ][line width=0.75]    (10.93,-3.29) .. controls (6.95,-1.4) and (3.31,-0.3) .. (0,0) .. controls (3.31,0.3) and (6.95,1.4) .. (10.93,3.29)   ;
				\draw    (222,140) -- (398,140) ;
				\draw [shift={(400,140)}, rotate = 180] [color={rgb, 255:red, 0; green, 0; blue, 0 }  ][line width=0.75]    (10.93,-3.29) .. controls (6.95,-1.4) and (3.31,-0.3) .. (0,0) .. controls (3.31,0.3) and (6.95,1.4) .. (10.93,3.29)   ;
				\draw [shift={(220,140)}, rotate = 0] [color={rgb, 255:red, 0; green, 0; blue, 0 }  ][line width=0.75]    (10.93,-3.29) .. controls (6.95,-1.4) and (3.31,-0.3) .. (0,0) .. controls (3.31,0.3) and (6.95,1.4) .. (10.93,3.29)   ;
				\draw    (340,100) -- (320,140) ;
				\draw [shift={(320,140)}, rotate = 116.57] [color={rgb, 255:red, 0; green, 0; blue, 0 }  ][fill={rgb, 255:red, 0; green, 0; blue, 0 }  ][line width=0.75]      (0, 0) circle [x radius= 3.35, y radius= 3.35]   ;
				\draw [shift={(340,100)}, rotate = 116.57] [color={rgb, 255:red, 0; green, 0; blue, 0 }  ][fill={rgb, 255:red, 0; green, 0; blue, 0 }  ][line width=0.75]      (0, 0) circle [x radius= 3.35, y radius= 3.35]   ;
				
				\draw (301,185) node [anchor=north west][inner sep=0.75pt]   [align=left] {$\displaystyle p_{2}$};
				\draw (352,83) node [anchor=north west][inner sep=0.75pt]   [align=left] {$\displaystyle p_{3}$};
				\draw (205,102) node [anchor=north west][inner sep=0.75pt]   [align=left] {$\displaystyle L=\ H_{1} \ \cap Q$};

			\end{tikzpicture}
		\caption*{$L$ is in $\Omega_{2,2}$}
		\end{minipage}
		\caption{Geometry in $Q \cong \bP^2$}
		\label{fig:twoschubert}
	\end{figure}
	
	Finally, we note that $\cV = \Omega_1^3 \subset \Gr(2,n)$ will have these recovered components for all $n \geq 5$, since nothing about the geometry in $Q$ uses the value of $n$.
	\end{exa}

    \begin{exa}[Codimension seven variety]\label{eg:seven}
		Similarly, $\Omega_1^7 \subset \Gr(2,n)$ will have recovered Schubert components of two classes: $\binom{7}{3}$ of type $\Omega_{6^2}$ (lines contained in $H_i \cap H_j \cap H_k$) and $\binom{7}{4}$ of type $\Omega_7$ (lines meeting $H_i \cap H_j \cap H_k \cap H_l$).

    To see this, we observe that the Chow-Lam locus lives in $\Gr(5,n).$ Consider a general subspace $Q$ in the Chow-Lam locus. Then $Q$ has dimension $4$ and contains seven planes $P_i:= Q \cap H_i.$ A generically chosen $Q$ will not contain a transversal to $P_1, ..., P_7.$ We would like to show that if $Q$ is \emph{specially} chosen to contain a line $L$ in one of the Schubert classes above, then $Q$ must also contain a transversal line $L'$ to the seven $P_i.$

    For an example of the first type, consider the case $n=8$ and let $L$ denote the line $H_1 \cap H_2 \cap H_3.$ Suppose $Q$ contains $L$. Then any line $L'$ in $Q$ which meets $L$ and $P_4, ..., P_7$ will be a common transversal to the seven planes. The subvariety of such $L'$ in $\Gr(2,5)$ has cohomology  class $[\Omega_2][\Omega_1]^4 \neq 0,$ and thus a common transversal must exist.
    
    For an example of the second type, consider the case $n=9$ and a line $L$ that meets the point $H_1 \cap ... \cap H_4.$ Suppose $Q$ contains $L,$ and thus the intersection of the four hyperplanes Then any line $L'$ in $Q$ through the point $H_1 \cap ... \cap H_4$ and meeting $P_5, P_6, P_7$ will be a common transversal. This has class $[\Omega_3][\Omega_1]^3,$ which is nonzero in $H^*(\Gr(2,5)).$ 
		
	Note that we have conditions on $n$ in both cases for these components to actually appear for generic $H_i$. For the class $\Omega_{6^2}$, we must have that the three $H_i$ intersect in at least a line in $\bP^{n-1}$, so that $L$ can be contained in the intersection. Thus we want $(n-1) - \text{codim } H_1 \cap H_2 \cap H_3  = (n-1) - 3 \cdot 2 \geq 1,$ meaning that $n$ is at least $8.$ For the second class, we want that $(n-1) - \text{codim } H_1 \cap H_2 \cap H_3 \cap H_4 \geq 0,$ so $n \geq 9.$ There is also a cohomological bound on $n$. We want the classes of $ 
		\Omega_{6^2}$ and $\Omega_7$ to be nonzero in $\Gr(2,n),$ so $n$ must be at least $8$ and $9$, respectively.
		
	\end{exa}

	It turns out that $\Omega_1^{ki+1} \subset \Gr(k,n)$ will have additional recovered Schubert components for all $i$, provided that $n$ is sufficiently high. Here we list some recovered components for small $k$ and $i$, assuming that $n$ is higher than some lower bound which we describe in Theorem \ref{thm:recovered}.

 \begin{table}[h!]
\centering
\begin{tabular}{||p{0.1cm} p{0.1cm} p{.6cm} p{2.1cm} p{4.2cm} p{1.8cm} ||} 
 \hline
 $k$ & $i$ & $\cV$ & Recovered & Interpretation, with $\cap_i H:=H_1\cap\ldots\cap H_i$ & Bounds on $n$ \\ [0.5ex] 
 \hline\hline
 2 & 3 & $\Omega_1^7$ & $\Omega_7, \Omega_{6^2}$ & line meets $\cap_4 H$, line contained in $\cap_3 H$ & $9, \ 8$ \\ 
 2 & 4 & $\Omega_1^9$ & $\Omega_9, \Omega_{8^2}$ & line meets $\cap_5 H$,  line contained in $\cap_4 H$ & $11, 10$ \\
 2 & i & $\Omega_1^{2i+1}$ & $\Omega_{2i+1}, \Omega_{(2i)^2}$ & line meets $\cap_{i+1}H$, line contained in $\cap_i H$& $2i+3,$ $2i+2$ \\
 3 & 4 & $\Omega_1^{13}$ & $\Omega_{13}, \Omega_{11^2}, \Omega_{9^3}$ & plane meets $\cap_5 H$, plane meets $\cap_4 H$ in a line, plane contained in $\cap_3 H$ & $16, 14, 12$\\
 3 & i & $\Omega_1^{3i+1}$ & $\Omega_{3i+1},$ $\Omega_{(3i-1)^2},$ $\Omega_{(3i-3)^3}$ & plane meets $\cap_{i+1} H$, plane meets $\cap_i H$ in a line, plane contained in $\cap_{i-1} H$ & $3i+4$, $3i+2$, $3i$  \\[1ex] 
 \hline
\end{tabular}
\caption{Some recovered components for some small values of $k$ and $i$}
\label{table:components}
\end{table}

\begin{thm}\label{thm:recovered}
		Fix $k$ and $i > k.$ Consider $V:= \Omega_1^{ki+1} \subset \Gr(k,n)$ for $n > k(i+1)+1$. Then the Chow-Lam recovery $W_V$ will contain recovered components of Schubert types $\Omega_{ki+1}, \Omega_{(k(i-1)+2)^2}, \Omega_{(k(i-2)+3)^3}, ..., \Omega_{(k(i-k+1)+k)^k}.$
	\end{thm}
	\begin{proof}
	Fix $a$ between $0$ and $k-1$, and fix any $i+1-a$ of the $ki+1$ codimension $k$ spaces $H_i$. Let $H$ denote their intersection of codimension $k(i+1-a)$ in $\bP^{n-1}$. Consider the Schubert variety of linear spaces of projective dimension $k-1$ which intersect $H$ in dimension at least $a.$ This is isomorphic to $\Omega_{(k(i-a)+a+1)^{a+1}}.$ We will show that this Schubert variety is a recovered component of $\Omega_1^{ki+1}.$
	
	For example, for $i=3$, $k=2$ and $a = 0,$ we may fix codimension $2$ spaces $H_1, H_2, H_3, H_4,$ and consider the Schubert variety of lines in $\bP^{n-1}$ which intersect the codimension $8$ intersection $H: = H_1 \cap H_2 \cap H_3 \cap H_4.$ 
 This is a Schubert variety isomorphic to $\Omega_7.$  It is nonempty only if $H$ is nonempty, i.e. $n\geq 9.$

  In fact, there are $(k-1)$-spaces that intersect $H$ in dimension at least $a$ only if the dimension of $H$ in $\bP^{n-1}$ is at least $a.$ Thus we must have 
    $$n-1 - \codim \ H \geq a,$$
    where $\text{codim } H = k(i+1-a).$ So $$n\geq k(i+1)+1-a(k-1).$$ This lower bound on $n$ is largest when $a=0$, in which case it is $k(i+1)+1$.
	
	Let us turn to the ambient Grassmannian of the Chow-Lam locus. The dimension of $\Omega_1^{ki+1}$ is $k(n-k) - (ki+1) = k(r-k)-1$ where $r = n-i.$ Thus its Chow-Lam locus is a subset of $\Gr(i+k, n).$ 
 
 Now, choose a $(k-1)$-space $P$ intersecting $H$ in at least dimension $a,$ so that $P$ is in the relevant Schubert variety. Next, choose a generic space $Q \cong \bP^{i+k-1}$ containing $P.$ Let $P_1, ..., P_{ki+1}$ be the intersections of the codimension $k$-spaces $H_1, ..., H_{ki+1}$ with $Q.$ We would like to show that the $P_i$ have a common transversal $P'\in \Gr(k,\hat Q)=\Gr(k,i+k)$ in $Q.$ 
	
	Without the assumption that $Q$ contains $P,$ the $ki+1$ linear spaces $P_i$ do not generically have such a transversal; indeed, meeting $P_i$ is a codimension one condition on $P' \in \Gr(k, i+k)$ and the total dimension of $\Gr(k, i+k)$ is $ki.$ However, if $Q$ contains $P,$ then $Q$ also intersects $H$ in at least dimension $a$ and any $(k-1)$-space $P'$ in $Q$ meeting $H$ and the remaining $P_i$ is a common transversal. The class  of such $P'$ is $\Omega_{i-a} \cdot \Omega_1^{ki + 1 - (i+1-a)},$ which is nonzero in the cohomology ring of $\Gr(k, i+k).$ Thus there exists a common transversal.

	\end{proof}

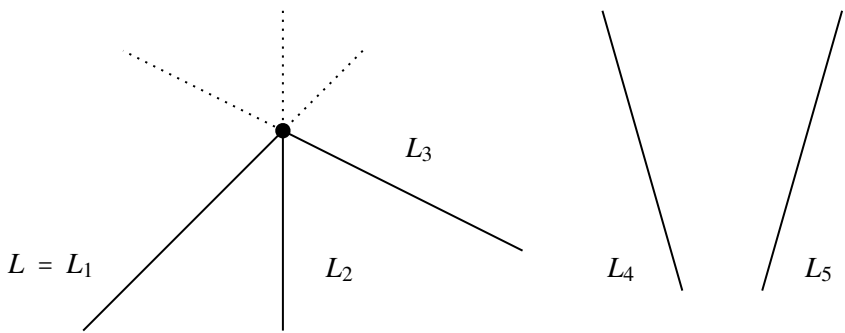
\begin{figure}[!h]
		\centering
		
		\tikzset{every picture/.style={line width=0.75pt}} 
		
		\begin{tikzpicture}[x=0.75pt,y=0.75pt,yscale=-1,xscale=1]
			
			\draw    (220,100) -- (120,200) ;
			\draw [shift={(220,100)}, rotate = 135] [color={rgb, 255:red, 0; green, 0; blue, 0 }  ][fill={rgb, 255:red, 0; green, 0; blue, 0 }  ][line width=0.75]      (0, 0) circle [x radius= 3.35, y radius= 3.35]   ;
			\draw    (220,100) -- (220,200) ;
			\draw    (340,160) -- (220,100) ;
			\draw  [dash pattern={on 0.84pt off 2.51pt}]  (220,100) -- (140,60) ;
			\draw  [dash pattern={on 0.84pt off 2.51pt}]  (260,60) -- (220,100) ;
			\draw  [dash pattern={on 0.84pt off 2.51pt}]  (220,40) -- (220,100) ;
			\draw    (380,40) -- (420,180) ;
			\draw    (500,40) -- (460,180) ;
			
			\draw (81,160) node [anchor=north west][inner sep=0.75pt]   [align=left] {$\displaystyle L\ =\ L_{1}$};
			\draw (240,162) node [anchor=north west][inner sep=0.75pt]   [align=left] {$\displaystyle L_{2}$};
			\draw (280,102) node [anchor=north west][inner sep=0.75pt]   [align=left] {$\displaystyle L_{3}$};
			\draw (381,162) node [anchor=north west][inner sep=0.75pt]   [align=left] {$\displaystyle L_{4}$};
			\draw (480,162) node [anchor=north west][inner sep=0.75pt]   [align=left] {$\displaystyle L_{5}$};

		\end{tikzpicture}
		\caption{Geometry of the lines $L_i$ in $Q \cong \bP^3$}
		\label{fig:schubert7}
	\end{figure}

	We caution that not all recovered components of $\Omega_1^{ki+1} \subset \Gr(k,n)$  are of the form in Theorem \ref{thm:recovered}. For example, the Schubert variety $\Omega_{2,2}\cap \Omega_3$ of lines contained in $H_1$ and intersecting $H_2 \cap H_3$ is in the recovered locus of $\Omega_1^5 \subset \Gr(2,7).$ Indeed, consider a line $L$ in this Schubert variety. In a generic $3$-space $Q \cong \bP^3$, there are five lines $L_i := H_i \cap Q.$ If $Q$ contains $L$, then $L_1 = L$ intersects $H_2 \cap H_3$ at a point $p.$ Then $p$ must be in $Q$, so it is also the intersection point of the lines $L_2$ and $L_3.$ Thus there are common transversals of the five lines: those which pass through the point $p$ and intersect the lines $L_4, L_5$, of which there are one since $[\Omega_2] \cdot [\Omega_1]^2=1$ in $\Gr(2,\hat Q)$. This is shown in Figure \ref{fig:schubert7}. However, this case is special because $i$ is small, which allows for the possibility that $L_i$ is equal to $L$; for $i$ large,  $H_i \cap Q$ will have dimension higher than $k.$

	\begin{conj}
		If $i \geq k$ then all recovered components of $\Omega_1^{ki}$ are of the form in Theorem~\ref{thm:recovered}.
	\end{conj}

	One could also ask what happens when a linear section of $\Gr(k,n)$ is not the intersection of Schubert hyperplane sections: do generic $(k(r-k)-1)$-dimensional linear sections of the Grassmannian have residual components in their Chow-Lam recovery? This is difficult to establish in general. The special case of curves may indicate what to expect.
  Let $\cV$ be the intersection of $\Gr(2,n)$ with a general linear space of codimension $2n-5$ in $\bP^{\binom{n}{2}-1}.$ Then $\cV$ is a curve, and if $n \neq 4$, there is no line $L$ in $X_\cV$ that does not belong to $\cV$, i.e. $W_\cV=\cV.$  This is because $X_\cV$ is smooth and therefore contains a line not in the ruling only if $\cV$ is rational.  But, by adjunction, the canonical divisors on $\cV$ are the restriction of hypersurfaces of degree $2n-5-n=n-5$, so $\cV$ is rational only if $n=4$. We note that in this case $\cV$ must have degree $2.$

 A threefold example shows that the analysis in general may be more involved.  In it we are able to prove there are no residual components, using the tools we developed in Section \ref{sec:criteria}.

 \begin{prop}
         Let $Q$ be the intersection of three generic hyperplanes in the Pl\"ucker space $\bP^9$ of $\Gr(2,5)$, and let $\cV = Q \cap \Gr(2,5).$ Then $W_\cV$ is equal to $\cV,$ i.e. has no residual components. 
     \end{prop}
     \begin{proof}
     Let $\cV$ be $Q \cap \Gr(2,5).$ Then $\cV$ is a three-fold and $r = 4$. Suppose there is a recovered line $L.$ By Theorem \ref{thm:resid}, we must have that the dimension of the space $\cV_0$ of lines in $\cV$ meeting $L$ is at least two. 
     
     Since $Q$ is generic, the variety $\cV$ is a smooth Fano threefold of degree $5$. In particular it has Picard rank one with Picard group generated by the class of a hyperplane section (cf. \cite[Corollary 6.6]{Isk}).
     Note that $\cV_0$ is a divisor in $\cV.$ Thus its class is some multiple of the hyperplane class in $\cV.$  

     We can describe $\cV_0$ even more explicitly. Let $T_L$ denote the tangent space $T_L \Gr(2,5)$ at the point $L.$ The variety $C := T_L \cap \Gr(2,5)$ is four-dimensional, and is a cone over $\bP^1 \times \bP^2 \subset \bP^5$ in $\bP^6.$  It is the four-dimensional Schubert variety $\Omega_3$ of lines in $\Gr(2,5)$ meeting $L.$ Since $\cV$ is smooth, the linear space $Q$ does not contain $[L]$.  So $Q\cap C$ is the intersection of $\bP^1 \times \bP^2 \subset \bP^5$ with at least one hyperplane. Since $\cV_0$ is a surface in this intersection, it has degree at most $3$.  But the degree of every surface in $\cV$ is divisible by five, so this is a contradiction.    
     \end{proof}

Based on these first examples we boldly conjecture:

     \begin{conj} Fix $r$, with $k<r\leq n$.
         Let  $\cV$ be the intersection of $\Gr(k,n) \subset \bP^N$ with a generic $(k(n-r)+1)$-dimensional linear space in $\bP^n.$ Then the algebraic sets $W_\cV$ and $\cV$ coincide.
     \end{conj}

     Finally, we return to the study of positroid varieties laid out in the introduction. These are boundary components of the positive Grassmannian. Their ideals are given explicitly by setting certain Pl\"ucker coordinates to zero. Understanding projections of these varieties would aid in understanding the amplituhedron itself. The class of examples in Theorem \ref{thm:recovered} includes certain positroid varieties, as some of them can be cut out by generic Schubert divisors.

     \begin{cor}
         Let $n = 2i+1$ be odd and at least three. Consider the positroid variety $\cV$ in $\Gr(2,2n)$ given by the vanishing of $p_{12}, ..., p_{2n-1, 2n}.$  Then $W_\cV \neq \cV.$
     \end{cor}

     \begin{proof}
         Each condition of the form $p_{i, i+1} = 0$ defines a Schubert hyperplane. Thus $\cV$ is of the form $\Omega_1^n,$ and we only need to check the bounds in Theorem \ref{thm:recovered}. But $2n > 2(i+1) + 1 = n + 2$ whenever $n$ is greater than two.
     \end{proof}

      It would be interesting to undertake a systematic study of $W_\cV$ where $\cV$ is a positroid variety, using the tools developed in this paper.


\begin{thebibliography}{10}
		\begin{small}
			\setlength{\itemsep}{-0.6mm}
   

            \bibitem{AHT} N.~Arkani-Hamed and J.~Trnka: {\em The Amplituhedron}, J.~High Energy Physics~{\bf 10} (2014) 30. 
			
			\bibitem{Catanese} F.~Catanese:
			{\em Chow Varieties, Hilbert Schemes, and Moduli Spaces of General Type},
			J.~Algebraic Geom.~{\bf 1} (1992) 561-595.

            \bibitem{Cayley} A.~Cayley: {\em On a new analytical representation of curves in space}, Quarterly Journal of Pure and Applied Mathematics (1860) 3:225-236.

            \bibitem{CW}
            W.-L.~Chow and B.~van der Waerden: {\em Zur algebraischen Geometrie. IX. \"Uber zugeordnete
            Formen und algebraische Systeme von algebraischen Mannigfaltigkeiten},
            Mathematische Annalen {\bf 113} (1937) 692--704.

            \bibitem{Fulton} W. Fulton and P. Pragacz: {\em Schubert Varieties and Degeneracy Loci}, Lecture Notes in Mathematics, Springer Berlin Heidelberg (1998).

            \bibitem{GKZ} I.M Gelfand and M.M Kapranov and A.V Zelevinsky: 
            {\em Hyperdeterminants}, Advances in Mathematics~{\bf 96} (1992) 2:226-263.


\bibitem{Isk} V.A. Iskovskih: 
            {\em Fano 3-folds I,}, Math. USSR Izv.~{\bf 11} (1977) 485-527.

            \bibitem{Lam} T. Lam: {\em Amplituhedron cells and Stanley symmetric functions}, Commun. Math. Phys.~{\bf 343} (2016) 1025–1037 
            
			\bibitem{PS} E.~Pratt and B.~Sturmfels:
			{\em The Chow-Lam Form},
			{\tt arXiv:2401.10795}.

    \bibitem{Williams} L.~Williams: {\em The positive Grassmannian, the amplituhedron, and cluster algebras}, International Congress of Mathematicians~{\bf 10} (2021) 

		\end{small}
	\end{thebibliography}
\end{document}